\documentclass{amsart}

\usepackage{amsmath, amssymb, amsthm, esint}

\newcommand{\bmop}{\text{BMO}_{\gamma}}
\newcommand{\bmonorm}[1]{\|#1\|_{\bmop}}

\DeclareMathOperator*{\esssup}{ess\,sup}

\newcommand{\hilb}{H_{\gamma}}
\newcommand{\comm}[2]{[#1, #2]}
\newcommand{\commbh}{\comm{b}{\hilb}}

\newcommand{\re}{\operatorname{Re}}
\newcommand{\im}{\operatorname{Im}}

\numberwithin{equation}{section}

\newtheorem{theorem}{Theorem}[section]
\newtheorem{defn}{Definition}[section]
\newtheorem{prop}{Proposition}[section]

\title{Commutators of Hilbert transforms along monomial curves}
\author[Tyler Bongers]{Tyler Bongers}
\address{Tyler Bongers, Department of Mathematics \& Statistics, Washington University in Saint Louis, One Brookings Drive Saint Louis, MO 63130--4899}
\email{tyler.bongers@wustl.edu}
\author[Zihua Guo]{Zihua Guo}
\address{Zihua Guo, School of Mathematical Sciences, Monash University, Melbourne, VIC 3800, Australia}
\email{zihua.guo@monash.edu}
\author[Ji Li]{Ji Li}
\address{Ji Li, Department of Mathematics, Macquarie University, NSW, 2109, Australia}
\email{ji.li@mq.edu.au}
\author[Brett D. Wick]{Brett D. Wick}
\address{Brett D. Wick, Department of Mathematics \& Statistics, Washington University in Saint Louis, One Brookings Drive Saint Louis, MO 63130--4899}
\email{wick@math.wustl.edu}
\date{\today}

\begin{document}

\begin{abstract}
The Hilbert transforms associated with monomial curves have a natural non-isotropic structure. We study the commutator of such Hilbert transforms and a symbol $b$ and prove the upper bound of this commutator when $b$ is in the corresponding non-isotropic  BMO space by using the Cauchy integral trick. We also 
consider the lower bound of this commutator by introducing a new testing BMO space associated with the 
given monomial curve, which shows that the classical non-isotropic  BMO space is contained in the testing BMO space. We also show that the non-zero curvature of such monomial curves are important, since when 
considering Hilbert transforms associated with lines, the parallel version of non-isotropic  BMO space and testing BMO space have overlaps but do not have containment.
\end{abstract}

\maketitle

%\tableofcontents

\section{Introduction}

Given an operator $T$ and a function $b$, the commutator $[b, T]$ is defined formally by the equality
$$[b, T] f = b T(f) - T(bf)$$
for all functions in an appropriate function space. If $T$ is an operator which is bounded between, e.g. $L^2$ and itself, then there is a natural bound
$$\|[b, T] : L^2 \to L^2\| \le 2 \|b\|_{\infty} \|T : L^2 \to L^2\|.$$
However, this bound is generally not optimal since it ignores any potential cancellation that may exist within the commutator. In the case that $T$ is a singular integral operator of the form $Tf(x) = \int_{\mathbb{R}^n} K(x, y) f(y) \, dy$, we have a representation
$$[b, T] f(x)  = \int_{\mathbb{R}^n} \left(b(x) - b(y)\right) K(x, y) f(y) \, dy.$$
There is typically cancellation carried by the kernel $K$, but now there is an additional cancellation due to the oscillation of $b$. As such, it is natural to search for bounds on commutators that involve the space BMO of functions with bounded mean oscillation; these are functions for which
$$\|b\|_{BMO} := \sup_Q \int_Q |b - b_Q| \, dx < \infty$$
for all cubes $Q$; here, $b_Q$ denotes the average of $b$ over $Q$.

In the classical setting, the seminal paper \cite{CoiRocWei76} of Coifman, Rochberg, and Weiss shows that whenever $b \in BMO$, the commutator of $b$ with any Riesz transform on $\mathbb{R}^n$ is bounded from $L^2(\mathbb{R}^n)$ to itself; recall that the Riesz transform $R_j$ is given by 
$$R_j f(x) = \text{p.v. } c_n \int_{\mathbb{R}^n} \frac{x_j - y_j}{|x - y|^{n + 1}} f(y) \, dy$$
for $j \in \{1, ..., n\}$. Moreover, there is a converse result: the Riesz transforms characterize BMO in the sense that a symbol $b$ such that $[b, R_j]$ is $L^2$ bounded for all $j$ must lie in BMO. Quantitatively, we have that (up to a purely dimensional constant)
$$\|b\|_{BMO} \sim \sum_{j = 1}^n \|[b, R_j] : L^2 \to L^2\|.$$

In this paper, we will consider an integral operator with a singularity that is not induced by a classical Calder\'on-Zygmund kernel. The \emph{parabolic Hilbert transform} of a function $f : \mathbb{R}^2 \to \mathbb{C}$ is defined by 
$$H_{\gamma} f(x) = \text{p.v.} \int_{-\infty}^{\infty} f(x - \gamma(t)) \, \frac{dt}{t}$$
where $\gamma(t) = (t, t^2)$. This operator arises naturally in the study of the role of curvature (see e.g. \cite{SteWai78} and the discussion in \cite[Section 8.8]{Duo01}) and parabolic differential operators. This operator is substantially different from a Calder\'on-Zygmund operator, since the singularity is localized to a curve. Although it is known that $H_{\gamma}$ is bounded as an operator on $L^p(\mathbb{R}^2)$ for all $p$ in the reflexive range (see, e.g. \cite[Section 8.7]{Duo01}), the geometry is rather subtle. For a particular example of the complications involved, note that the parabolic Hilbert transform of an indicator function depends on the shape of the underlying set; this is not the case in the classical setting \cite{Lae12}. Moreover, classical techniques such as the median method (see, e.g. \cite{Hyt18} for the usage of this) break down due to the complicated geometry of the operator. Ultimately, this complicates the search for a BMO-type characterization for bounded commutators.

A key tool to study singular integral operators is the idea of sparse domination; recently, there has been much focus on refining classical boundedness results to sparse domination results. In this case, it is necessary to adapt the idea to the geometry of the operator. We will call a set $Q$ a parabolic cube if $Q = I \times J$ for intervals $I$ and $J$ with $|J| = |I|^2$; we write $\ell(Q) = |I|$. A collection $\mathcal{P}$ of parabolic cubes is called $\delta$-sparse (with $\delta \in (0, 1)$) if there exist pairwise disjoint sets $E_P \subset P$ with $|E_P| \ge \delta |P|$ for all $P \in \mathcal{P}$. Correspondingly, a sparse form $\Lambda_{\mathcal{P}, r, s}$ is defined by
$$\Lambda_{\mathcal{P}, r, s}(f, g) = \sum_{P \in \mathcal{P}} |P| \left(\frac{1}{|P|} \int_P |f|^r\right)^{1/r} \left(\frac{1}{|P|} \int_P |g|^s\right)^{1/s}$$
for compactly supported smooth functions $f$ and $g$. The sparse domination theorem of Cladek and Ou \cite{ClaOu18} shows that the parabolic Hilbert transform is dominated by a sparse form for a wide range of exponents. This sparse domination result then gives a direct relationship between the operator and a family of cubes; once an appropriate sparse form is involved, there are known weighted bounds that can be adapted to the commutator setting. This suggests that the appropriate BMO space necessary for boundedness of commutators $[b, H_{\gamma}]$ should be adapted to collection of parabolic cubes; it turns out that this is indeed the case.

Following this idea, we will say that $b$ is in parabolic BMO if 
$$\|b\|_{BMO_{\gamma}} := \sup_{P \in \mathcal{P}} \frac{1}{|P|} \int_P |b - b_P| \, dx < \infty$$
where $\mathcal{P}$ is the collection of all parabolic cubes.   The main result of this paper is the following theorem:
\begin{theorem}\label{theorem:parabolic_sufficient}
If $b \in BMO_{\gamma}$, and $1<p<\infty$ then the commutator $[b, H_{\gamma}]$ is bounded from $L^p(\mathbb{R}^2)$ to itself, and
$$\|[b, H_{\gamma}] : L^p(\mathbb{R}^2) \to L^p(\mathbb{R}^2)\| \lesssim \|b\|_{BMO_{\gamma}}.$$
\end{theorem}
Additionally, we will give a partial converse result based on a similar BMO-type space. Given a parabolic cube $P$, consider a set $E_P = \{x - \gamma(t) : x \in P, 9 \ell(P) \le t \le 10 \ell(P)\}$ and the interval $I_{x, E_P} = \{t : x - \gamma(t) \in E_P\}$. If $\mu$ is the Haar measure $dt/t$ on $\mathbb{R}^+$, the testing-BMO norm is defined to be
$$\|b\|_{\text{test}} = \sup_{P \in \mathcal{P}} \frac{1}{|P|} \int_P \left|b(x) - \frac{1}{\mu(I_{x, E_P})} \int_{I_{x, E_P}} b(x - \gamma(t)) \, d\mu(t)\right| \, dx<\infty.$$
This captures a notion of oscillation of $b$ along the parabolic curves of the operator $H_{\gamma}$. We then have
\begin{theorem}\label{theorem:parabolic_necessary}
If the commutator $[b, H_{\gamma}]$ is bounded from $L^2(\mathbb{R}^2)$ to itself, $b$ is a testing-BMO function with
$$\|b\|_{\text{test}} \lesssim \|[b, H_{\gamma}] : L^2(\mathbb{R}^2) \to L^2(\mathbb{R}^2)\|.$$
\end{theorem}

The following is an outline of the paper. In Section 2, we will use the Cauchy integral trick to establish the boundedness of the commutator $[b, H_{\gamma}]$. A key tool is the sparse domination result of Cladek and Ou, together with weighted bounds for sparse operators due to Bernicot, Frey, and Petermichl \cite{BerFrePet16}. In Section 3, we will give the partial converse result involving testing BMO; this will involve as close an adaptation of the median method as the geometry of this operator allows. We also discuss the obstacles in improving this to a full converse. In Section 4, we extend the results of the previous sections to related operators with the same basic geometric structure; in particular, our techniques apply equally well to monomial curves in $\mathbb{R}^n$ and torsion-free curves. Finally, Section 5 considers the necessary role that curvature plays in these estimates.

\section{Sufficient condition for boundedness: Proof of Theorem \ref{theorem:parabolic_sufficient}}
In this section, we will show that if $b \in \bmop$, then there is a universal constant $C$ depending only on $p$ such that 
$$\|[b, H_{\gamma}] : L^p(\mathbb{R}^2) \to L^p(\mathbb{R}^2)\| \le C \bmonorm{b}.$$  
Hence forth we write $L^p:=L^p(\mathbb{R}^2)$ to shorten notation.

\begin{proof}
We will begin with $p = 2$. Since both sides of the claimed inequality scale linearly in the $\bmop$ norm of $b$, we assume without loss of generality that $\bmonorm{b} = 1$. The proof will use the Cauchy integral trick originally observed by Coifman, Rochberg, and Weiss \cite{CoiRocWei76}. Fix $z \in \mathbb{C}$ and form the operator $e^{zb/2} H_{\gamma} e^{-zb/2}$; expanding in a formal power series, we have
\begin{align*}
e^{zb/2} H_{\gamma}(e^{-zb/2}f) &= \left(1 + \frac{zb}{2} + \cdots\right) H_{\gamma} \left(\left(1 - \frac{zb}{2} + \cdots\right) f\right) \\
&= H_{\gamma} f + \frac{zb}{2} H_{\gamma} f + H_{\gamma}\left(-\frac{zb}{2} f\right) + z^2 T \\
&= H_{\gamma} f + \frac z 2 [b, H_{\gamma}] f + z^2 T
\end{align*}
for some operator $T$. Differentiating in $z$, we therefore have
$$\commbh f = 2 \left. \frac{d}{dz}\right|_{z = 0} e^{zb/2} \hilb(e^{-zb/2} f).$$
Fix a positive quantity $\epsilon$, to be determined later. The Cauchy integral formula gives the representation
$$
\commbh f = \frac{1}{\pi i} \int_{\partial \mathbb{D}_{\epsilon}} \frac{e^{zb/2} \hilb(e^{-zb/2} f)}{z^2} \, dz
$$
over the boundary of the disk $\mathbb{D}_{\epsilon}$. We now take $L^2$ norms and estimate the integral roughly:
\begin{align}
\|\commbh f\|_{L^2} &\le \frac 1 {\pi} \int_{\partial \mathbb{D}_{\epsilon}} \frac{\|e^{zb/2} H_{\gamma}(e^{-zb/2} f)\|_{L^2}}{|z|^2} \, |dz| \nonumber \\
&\lesssim \frac{1}{\epsilon} \sup_{|\re z| = \epsilon} \|e^{zb/2} H_{\gamma}(e^{-zb/2} f)\|_{L^2}. \label{eqn:weighted_hilbert_bound}
\end{align}
This allows a reduction to the study of weighted estimates. In particular, if we define $\tilde{f} = e^{-b \im z / 2} f$,
\begin{align*}
\|e^{zb/2} H_{\gamma}(e^{-zb/2} f)\|_{L^2} &= \int_{\mathbb{R}^2} |e^{zb/2}|^2 |H_{\gamma}(e^{-zb/2} f)|^2 \, dx \\
&= \int_{\mathbb{R}^2} e^{b \re z} |H_{\gamma}(e^{-b \re z} \tilde{f})|^2 dx.
\end{align*}
Note that $|\tilde{f}| = |f|$ at all points, since we have only modified $f$ by a rotation. We now can define a positive weight $w = e^{b \re z}$ depending on $\re z$ and recognize that
\begin{align*}
\|e^{zb/2} H_{\gamma}(e^{-zb/2} f)\|_{L^2} &= \|H_{\gamma}(w^{-1/2} \tilde{f})\|_{L^2(w)} \\
&\le \|H_{\gamma} : L^2(w) \to L^2(w)\| \cdot \|w^{-1/2} \tilde{f}\|_{L^2(w)} \\
&= \|H_{\gamma} : L^2 (w) \to L^2(w)\| \cdot \|f\|_{L^2}
\end{align*}
Combining this estimate with \eqref{eqn:weighted_hilbert_bound}, we have for each $\epsilon > 0$ that
\begin{equation}\label{eqn:weighted_hilbert}
\|[b, H_{\gamma}] f\|_{L^2} \lesssim \left( \sup_{|\re z| = \epsilon} \frac{\|H_{\gamma} : L^2(w) \to L^2(w)\|}{\epsilon} \right) \|f\|_{L^2}
\end{equation}
Therefore, it is necessary to find a weighted bound for the Hilbert transform. Typically, weighted bounds involve the $A_p$ characteristic of the weight, which in turn is related to the BMO norm of $\log w$. Finally, since $w = e^{b \re z}$, this will give the desired estimate. There are several complications, however; since $H_{\gamma}$ is a more complicated operator than a classical Calder\'on-Zygmund operator, the corresponding weighted bounds are more involved (and the $A_p$ characteristic will involve parabolic cubes rather than the typical dyadic grid). In this case, the relevant estimates here pass through the sparse domination results of Cladek and Ou \cite{ClaOu18}, and are originally due to Bernicot, Frey and Petermichl \cite{BerFrePet16}. Following \cite[Corollary 1, Section 4]{ClaOu18}, fix exponents $r$ and $s$ such that the point $(1/r, 1/s')$ lies in the interior of the triangle with vertices 
$$(0, 0), (1, 0), (2/3, 1/3).$$
If $\alpha = \max\{\frac{1}{2 - r}, \frac{s' - 1}{s' - 2}\}$, the following estimate holds in terms of the Muckenhoupt and reverse H\"older characteristics of $w$:
\begin{equation}\label{eqn:cladek_ou}
\|\hilb : L^2(w) \to L^2(w)\| \le C(r, s) \left([w]_{A_{2/r}} [w]_{RH_{(s'/2)'}}\right)^{\alpha}.
\end{equation}

The exact selection of the parameters $r$ and $s$ will be made later. The next step is to compute the characteristics of $w$ in terms of the BMO norm. This has been done for BMO functions where the norm is defined relative to the typical dyadic grid; we remark here that there are no issues in carrying the classical arguments out in the parabolic setting. In particular, the grid of parabolic cubes still has a differentiation property (that is, parabolic averages of a continuous function converge to the value of the function as the scale tends to zero) and it is still possible to run stopping time arguments. 

Therefore, we need to estimate both the $A_{2/r}$ characteristic and the reverse H\"older characteristic of $w = e^{b \re z}$. There is a classical correspondence between $A_2$ and $\bmop$ via exponentiation; for a quantitative version of this, see B\'enyi, Martell, Moen, Stachura, and Torres \cite[Lemma 3.5]{BenMarMoeStaTor17}. If $g \in \bmop$, $p \in [1, \infty)$, and $|\lambda| \le \min\{1, p - 1\} \|g\|_{\bmop}^{-1}$ we have the estimate
$$[e^{\lambda g}]_{A_p} \le 4^{|\lambda| \|g\|_{\bmop}}.$$
Apply this to the function $b$ (with $\|b\|_{\bmop} = 1$) at exponent $p = 2/r$. Provided that $|\re z| \le 2/r - 1 < 1$, we then have the estimate
\begin{equation}\label{eqn:ap_estimate}
[w]_{A_{2/r}} = [e^{xb}]_{A_{2/r}} \le 4^{(2/r - 1) \cdot 1} < 4.
\end{equation}
Note further that only using $|\re z| \le 1$ gives the weaker estimate
$$[w]_{A_2} = [e^{b\re z}]_{A_2} \le 4.$$

We now turn to the reverse H\"older characteristic. The quantitative result that we need is by Hyt\"onen, Perez, and Rela \cite[Theorem 2.3]{HytPerRel12}: if $\eta \in A_{\infty}$, 
$$[\eta]_{RH_{1 + \sigma}} \le 2$$
for all $\sigma \le C [\eta]_{A_{\infty}}^{-1}$ for some positive constant $C$. (As remarked above, this theorem is proved only for the standard setting of dyadic cubes; however, the proof immediately translates into the parabolic setting with only a change of constants. Alternatively, we can apply \cite[Theorem 1.2]{HytPerRel12} which proves a similar result in the setting of a space of homogeneous type.) In any case, applying this to the weight $w \in A_2$ and noting that $[w]_{A_{\infty}} \le [w]_{A_2} \le 4$ for all $|x| \le 1$, we have that
\begin{equation}\label{eqn:rh_estimate}
[w]_{RH_{1 + \sigma}} \le 2
\end{equation}
for all $\sigma$ sufficiently small.

We are now ready to select the parameters. To choose $(r, s)$ in the acceptable parameter range of Cladek and Ou from equation \eqref{eqn:cladek_ou}, begin with $s$ sufficiently close to $1$ that $(s'/2)' < 1 + \sigma$; then select $r$ so that $(1/r, 1/s')$ is strictly within the interior of the parameter range. Note that both of these choices are independent of all other parameters, such as $b$. If $\epsilon := 2/r - 1 > 0$, then the $A_{2/r}$ estimate can be carried out as in \eqref{eqn:ap_estimate}. Combining the weighted bound \eqref{eqn:weighted_hilbert}, the mixed Muckenhoupt-reverse H\"older estimate \eqref{eqn:cladek_ou}, and the characteristic estimates \eqref{eqn:ap_estimate} and \eqref{eqn:rh_estimate}, we have
\begin{align*}
\|[b, H_{\gamma}] : L^2 \to L^2\| &\lesssim \sup_{|\re z| \le \epsilon} \frac{\|H_{\gamma} : L^2(w) \to L^2(w)\|}{\epsilon} \\
&\lesssim_{r, s} \sup_{|\re z| \le \epsilon} \left([w]_{A_{2/r}} [w]_{RH_{(s'/2)'}}\right)^{\max\left\{\frac{1}{2 - r}, \frac{s' - 1}{s' - 2}\right\}} \\
&\le 8^{\max\left\{\frac{1}{2 - r}, \frac{s' - 1}{s' - 2}\right\}} \\
&\lesssim_{r, s} 1
\end{align*}
as desired and the proof is complete for $p = 2$.

We now turn to the case of $p \in (1, \infty)$. Again normalize $\bmonorm{b} = 1$. Proceeding as before with the Cauchy integral trick, we have the starting estimate
$$\|[b, H_{\gamma}] : L^p \to L^p\| \lesssim \sup_{|\re z| \le \epsilon} \frac{\|H_{\gamma} : L^p(w) \to L^p(w)\|}{\epsilon}$$
with $w = e^{xb}$. If $p \in (r, s')$, Corollary 1 of \cite{ClaOu18} then gives
$$\|[b, H_{\gamma}] : L^p \to L^p\| \lesssim_{r, s, p} \left([w]_{A_{p/r}} [w]_{A_{(s'/p)'}}\right)^{\alpha}$$
with $\alpha = \max\{\frac{1}{p - r}, \frac{s' - 1}{s' - p}\}$. Provided that we can find $\epsilon > 0$ depending only on $p$ and then bound $[w]_{A_{p/r}}$ and $[w]_{RH_{(s'/p)'}}$ universally as before, this will yield $L^p$ bounds.

For the Muckenhoupt characteristic, we have the estimate
$$[w]_{A_{p/r}} = [e^{b\re z}]_{A_{p/r}} \le 4^{\re z} \le 4$$
provided that $|\re z| \le \min\{1, p/r - 1\}$ (which will lead to the choice $\epsilon = p/r - 1$, depending only on the allowable parameters). For the reverse H\"older characteristic, note that $[w]_{A_{\infty}} \le 4$ and so for all $\sigma$ sufficiently small we again have the estimate
\begin{equation}\label{eqn:ap_other}
[w]_{RH_{1 + \sigma}} \le 2.
\end{equation}
Finally, note that the exponent $\alpha$ depends only on the allowable parameters.

We can now select the parameters. Choose $s$ very close to $1$ such that $(s'/p)'$ is sufficiently close to $1$ that the reverse H\"older estimate \eqref{eqn:ap_other} above holds; this selection depends only on $p$; note that $s'$ will be larger than $p$. Next choose $r$ such that $(1/r, 1/s')$ lies in the interior of the triangle with vertices at $(0, 0), (1, 1)$, and $(2/3, 1/3)$; as $s' > p$, this requires $r < p$. If $r$ and $s'$ are both sufficiently close to $p$, the point $(1/r, 1/s')$ will lie near the upper boundary of the triangle, but still within it. Finally, select $\epsilon = p/r - 1 > 0$. Since $r$ depends only on $s$, which depends only on $p$, all our parameters only depend on $p$. Thus,
\begin{equation}\label{eqn:lp_commutator}
\|[b, H_{\gamma}] : L^p \to L^p\| \lesssim_{p} \left([w]_{A_{p/r}} [w]_{A_{(s'/p)'}}\right)^{\alpha} \lesssim_{p} 1
\end{equation}
as desired.
\end{proof}

As an extension of these techniques, it is possible to extend the results to higher order commutators. The $k$-th order commutator $T^k_b$ is defined recursively by $T^0_b = H_{\gamma}$ and $T^{k + 1}_b = [b, T^k_b]$; note that $k = 1$ corresponds to the usual commutator $[b, H_{\gamma}]$. We then have the following boundedness result:
\begin{theorem}
If $b \in BMO_{\gamma}$, the $k$-th order commutator $T^k_b$ is bounded as an operator from $L^p(\mathbb{R}^2)$ to itself for all $p \in (1, \infty)$, and there is a constant $C = C(p)$ for which
$$\|T^k_b : L^p(\mathbb{R}^2) \to L^p(\mathbb{R}^2)\| \le (C^k\cdot k!) \|b\|_{BMO_{\gamma}}.$$
\end{theorem}

\begin{proof}
One can directly verify from the series representation that
$$T^k_b f = 2^k \left.\frac{d^k}{dz^k}\right|_{z = 0} e^{zb/2} H_{\gamma}(e^{-zb/2} f)$$
and therefore the Cauchy integral representation gives us
$$T^k_b f = \frac{2^k \cdot k!}{2\pi i} \int_{\partial \mathbb{D}_{\epsilon}} \frac{e^{zb/2} H_{\gamma}(e^{-zb/2} f)}{z^{k + 1}} \, dz.$$
Proceeding as in the first-order case, we have
$$\|T^k_b : L^p \to L^p\| \le 2^k \cdot k! \, \sup_{|x| \le \epsilon} \frac{\|H_{\gamma} : L^p(e^{xb/2}) \to L^p(e^{xb/2})\|}{\epsilon^k}.$$
We have already computed the norm of the weighted Hilbert transform in terms of the Muckenhoupt and reverse H\"older characteristics of the weight $e^{xb}$, and have shown that it can be bounded if $\epsilon$ is chosen to be sufficiently small. As stated in the leadup to equation \eqref{eqn:lp_commutator}, these bounds depend only on $p$; the claim follows.
\end{proof}

\section{Necessary condition for boundedness: Proof of Theorem \ref{theorem:parabolic_necessary}}
We will now give a BMO-type condition for $b$ which is necessary for the $L^2$ bounds of the commutator $[b, H_{\gamma}]$. Recall that for a parabolic cube $Q$ with dimensions $\ell(Q) \times \ell(Q)^2$, we have defined a set $E_Q$ by
$$E_Q = \{x - \gamma(t) : x \in Q, 9 \ell(Q) \le t \le 10 \ell(Q)\}.$$
For each $x \in Q$, let $I_{x, E_Q} = \{t : x - \gamma(t) \in E_Q\}.$ The set $E_Q$ represents the flow of the cube $Q$ along a parabola (on a time scale comparable to the length of $Q$), while $I_{x, E_Q}$ represents the time that the parabolic trajectory based at $x$ will spend in the flow $E_Q$. Our oscillation condition on $b$ is then
$$\|b\|_{\text{test}} := \sup_{Q \in \mathcal{P}}  \fint \left|b(x) - \frac{1}{\int_{I_{x, E_Q}} \frac{dt}{t}} \int_{I_{x, E_Q}} b(x - \gamma(t))\, \frac{dt}{t}\right| \, dx < \infty.$$
Alternatively, if $d\mu = dt / t$ is the Haar measure on $\mathbb{R}^+$, then
$$\|b\|_{\text{test}} := \sup_{Q \in \mathcal{P}} \fint \left|b(x) - \fint_{I_{x, E_Q}} b(x - \gamma(t)) \, d\mu(t)\right| \, dx < \infty.$$

We now turn to the proof of Theorem \ref{theorem:parabolic_necessary}; that is, that if $b$ is a locally integrable symbol, then
$$\|b\|_{\text{test}} \lesssim \|[b, H_{\gamma}] : L^2 \to L^2\|.$$

\begin{proof}
Suppose that $b$ is a locally integrable symbol for which the commutator $[b, H_{\gamma}]$ is bounded on $L^2$. Let $Q$ and $E$ be any two sets of positive and finite area; for any $x \in Q$, define the set
$$I_{x, E} = \{t \in \mathbb{R} : x - \gamma(t) \in E\}.$$
By the Cauchy-Schwarz inequality,
\begin{align}
\fint_Q |[b, H_{\gamma}] \chi_E| \, dx &\le \left(\fint_Q |[b, H_{\gamma}] \chi_E|^2 \, dx\right)^{1/2} \nonumber \\
&\le \frac{1}{|Q|^{1/2}} \left(\int_{\mathbb{R}^2} |[b, H_{\gamma}] \chi_E|^2 \, dx\right)^{1/2} \nonumber \\
&\le \frac{1}{|Q|^{1/2}} \|[b, H_{\gamma}] : L^2 \to L^2\| \cdot \|\chi_E\|_{L^2} \nonumber \\
&= \frac{|E|^{1/2}}{|Q|^{1/2}} \|[b, H_{\gamma}] : L^2 \to L^2\|. \label{eqn:testing_indicator}
\end{align}

On the other hand, we can compute the commutator in terms of the time sets $I_{x, E}$: 
\begin{align}
&\fint_Q |[b, H_{\gamma}] \chi_E| \, dx \nonumber \\
&= \fint_Q \left|b(x) \left(H_{\gamma} \chi_E\right)(x) - \left(H_{\gamma}(b \chi_E)\right)(x)\right| \, dx \nonumber \\
&= \fint_Q \left| \int_{-\infty}^{\infty}\left[ b(x) \chi_E(x - \gamma(t)) - b(x - \gamma(t)) \chi_E(x - \gamma(t)) \right]\, \frac{dt}{t} \right| \, dx \nonumber \\
&= \fint_Q \left| \int_{I_{x, E}} \left[b(x) - b(x - \gamma(t))\right] \, \frac{dt}{t} \right| \, dx \nonumber \\
&= \fint_Q \left( \int_{I_{x, E}} \frac{dt}{t} \right) \left|b(x) - \frac{1}{\int_{I_{x, E}} \frac{dt}{t}} \int_{I_{x, E}} b(x - \gamma(t)) \, \frac{dt}{t}\right| \, dx. \label{eqn:haar_oscillation}
\end{align}
Combining the inequalities \eqref{eqn:testing_indicator}--\eqref{eqn:haar_oscillation}, this yields
\begin{align*}
&\fint_Q \left|b(x) - \frac{1}{\int_{I_{x, E}} \frac{dt}{t}} \int_{I_{x, E}} b(x - \gamma(t)) \, \frac{dt}{t}\right| \, dx \\
&\le \frac{|E|^{1/2}}{|Q|^{1/2}} \left( \min_{x \in Q} \int_{I_{x, E}} \frac{dt}{t} \right)^{-1} \|[b, H_{\gamma}] : L^2 \to L^2\|.
\end{align*}
This chain of reasoning holds for all sets $Q$ and $E$ with positive and finite area. We now specialize: let $Q$ be a parabolic cube and $E = E_Q$. We will estimate the quantity $|E_Q|^{1/2} |Q|^{-1/2} \left(\min_{x \in Q} \int_{I_{x, E_Q}} dt/t\right)^{-1}$ and show that it is uniformly bounded away from $0$ and $\infty$; once this is established, it will yield a necessary condition for the symbol $b$ to have a bounded commutator. To this end, note that if $x \in Q$, then 
$$[9 \ell(Q), 10 \ell(Q)] \subseteq I_{x, E_Q} \subseteq [8 \ell(Q), 11\ell(Q)]$$
so that
$$\ln \frac{10}{9} \le \int_{I_{x, E_Q}} \frac{dt}{t} \le \ln \frac{11}{8}.$$
Moreover, by considering the flows of the corners of the parabolic cube $Q$, we can explicitly compute $E_Q$: it is a (non-parabolic) rectangle whose base length is $(10 - 9 + 1) \ell(Q) = 2 \ell(Q)$ and whose height is $(10^2 - 9^2 + 1) \ell(Q) = 20 \ell(Q)$; therefore, 
\begin{equation}\label{eqn:similar_measures}
|E_Q| \sim |Q|\end{equation}
holds with uniform constant. 
\end{proof}
Combining this theorem with the previous section, we have the inclusions
$$\text{parabolic BMO} \implies \text{bounded commutator} \implies \text{commutator testing}$$
and associated quantitative estimates
$$\|b\|_{\text{test}} \lesssim \|[b, H_{\gamma}] : L^2 \to L^2\| \lesssim \|b\|_{BMO_{\gamma}}.$$
However, it is unclear to us whether any of these inclusions are proper; as such, a full characterization of the symbols $b$ for which the commutator $[b, H_{\gamma}]$ is bounded has not yet been found. Note that it is possible to go directly from parabolic BMO to the testing condition without the intermediate commutator step:

\begin{prop} If $b \in BMO_{\gamma}$ then $b$ satisfies the commutator testing condition and
$$\|b\|_{\text{test}} \lesssim \|b\|_{BMO_{\gamma}}.$$
\end{prop}

\begin{proof}
Suppose that $b$ is a parabolic BMO function and fix $Q \in \mathcal{P}$. Let $Q_t := Q - \gamma(t)$ denote the parabolic flow of $Q$ by $t$, which is an offset parabolic cube from $Q$. We then have the estimate
\begin{align}
&\fint_Q \left|b(x) - \fint_{I_{x, E_Q}} b(x - \gamma(t)) \, d\mu(t)\right| \, dx \nonumber \\
&\le \fint_Q |b - b_Q| \, dx + \fint_Q \left|b_Q - \fint_{I_x, E_Q} b_{Q_t} \, d\mu(t)\right| \, dx \,+ \nonumber \\
&\quad\quad+ \fint_Q \left|\fint_{I_{x, E_Q}} b_{Q_t} \, d\mu(t) - \fint_{I_{x, E_Q}} b(x - \gamma(t)) \, d\mu(t)\right| \, dx \nonumber \\
&= \fint_Q |b - b_Q| \, dx + \fint_Q \left|\fint_{I_{x, E_Q}} b_Q - b_{Q_t} \, d\mu(t)\right| \, dx \,+ \nonumber \\
&\quad\quad+ \fint_Q \left|\fint_{I_{x, E_Q}} b_{Q_t} - b(x - \gamma(t)) \, d\mu(t)\right| \, dx \nonumber \\
&= (I) + (II) + (III). \label{eqn:bmo_sum_estimate}
\end{align}
The first term is controlled by the parabolic BMO norm of $b$. The second term of \eqref{eqn:bmo_sum_estimate} can be estimated by
$$\fint_Q \fint_{I_{x, E_Q}} |b_Q - b_{Q_t}| \, d\mu(t) \, dx \lesssim \fint_Q \fint_{I_{x, E_Q}} \|b\|_{BMO_{\gamma}} \, d\mu(t) \, dx = \|b\|_{BMO_{\gamma}}.$$
This is a consequence of the selection of $Q$ and $Q_t$ as equally sized parabolic cubes whose distance is bounded by a multiple of their (common) length. The averages of BMO functions on nearby cubes of comparable size can differ by at most a multiple of the BMO norm. Finally, the last term of \eqref{eqn:bmo_sum_estimate} can be handled by changing the order of integration and recognizing an average over a different parabolic cube:
\begin{align*}
&\fint_Q \left|\fint_{I_{x, E_Q}} b_{Q_t} - b(x - \gamma(t)) \, d\mu(t)\right| \, dx \\
&\le \fint_Q \fint_{I_{x, E_Q}} |b_{Q_t} - b(x - \gamma(t))| \, d\mu(t) \, dx \\
&\le \fint_Q \frac{1}{\int_{I_{x, E_Q}} \frac{dt}{t}} \int_{8 \ell(Q)}^{11\ell(Q)} |b_{Q_t} - b(x - \gamma(t))| \, d\mu(t) \, dx \\
&\le \frac{1}{\ln(10/9)} \int_{8\ell(Q)}^{11\ell(Q)} \fint_Q |b_{Q_t} - b(x - \gamma(t)| \, dx \, d\mu(t) \\
&\le \frac{\ln(11/8)}{\ln(10/9)} \sup_{8 \ell(Q) \le t \le 11 \ell(Q)} \fint_Q |b_{Q_t} - b(x - \gamma(t))| \, dx \\
&= \frac{\ln(11/8)}{\ln(10/9)} \sup_{8 \ell(Q) \le t \le 11 \ell(Q)} \fint_{Q_t} |b_{Q_t} - b| \, dx \\
&\lesssim \|b\|_{BMO_{\gamma}}
\end{align*}
as desired.
\end{proof}

As a further exploration of space of functions satisfying the testing inequality, note that we have an inequality reminiscent of a simple John-Nirenberg inequality in the special case of a simple symbol $b$.
\begin{prop}
Let $b = \sum_Q b_Q \chi_Q$ be a sum over a disjoint collection of cubes $\{Q\}$. If $Q$ and $R$ are adjacent parabolic cubes, then $|b_Q - b_R| \le \|b\|_{\text{test}}$.
\end{prop}

\begin{proof}
Without loss of generality, $R$ lies below or to the left of $Q$. It is possible to find a small parabolic cube $\tilde{Q} \subseteq Q$ such that $E_{\tilde{Q}} \subseteq R$. Since $b$ is constant on $Q$ and $R$, the testing condition applied to $\tilde{Q}$ reduces to the desired inequality.
\end{proof}

\section{The non-parabolic setting}

The techniques used in Sections 2 and 3 are not specialized to the parabolic setting; rather, they can be extended to a more general class of smooth curves. We will call a function $\eta: \mathbb{R} \to \mathbb{R}^n$ a monomial curve if $\eta$ is given by
$$\eta(t) = \left\{\begin{array}{cc} (\epsilon_1 |t|^{\alpha_1}, ..., \epsilon_n |t|^{\alpha_n}) & \text{ if } t \ge 0 \\
(\epsilon_1' |t|^{\alpha_1}, ..., \epsilon_n' |t|^{\alpha_n}) & \text{ if } t < 0
\end{array}\right\}$$
where the coefficients $\epsilon_i$ and $\epsilon_i'$ are chosen from $\{-1, 1\}$, and there is an index $j$ with $\epsilon_j \ne \epsilon_j'$. The associated Hilbert transform along $\eta$ is then
$$H_{\gamma}f(x) = p.v. \int_{\mathbb{R}} f(x - \gamma(t)) \, \frac{dt}{t}.$$

A monomial curve induces a family of so-called $\eta$-cubes and an associated BMO space. We will call a cube an $\eta$-cube if it is of the form $Q = I_1 \times \cdots I_n \subseteq \mathbb{R}^n$ and the side-lengths are related by
$$|I_1|^{\alpha_1} = \cdots = |I_n|^{\alpha_n}.$$
Naturally, there is a BMO space normed by
$$\|b\|_{BMO_{\eta}} = \sup_{Q \text{ an } \eta-\text{cube}} \frac{1}{|Q|} \int_Q |b - b_Q| \, dx.$$
As proven in \cite{ClaOu18}, there is also an associated sparse-domination theorem that directly generalizes the parabolic setting. This family of cubes has the same differentiation properties that the dyadic (and parabolic) cubes enjoy; moreover, the sparse domination theorem for $H_{\eta}$ involves constants only depending on the dimension $n$ and the curve $\eta$; the range of exponents where sparse domination holds does not depend on the curve $\eta$, but is purely dimensional. Therefore, we may apply the previous arguments \emph{mutatis mutandis} to prove commutator bounds for $H_{\eta}$ in the monomial context as well:
\begin{theorem}
Let $p \in (1, \infty)$. If $\eta$ is a monomial curve and $b \in BMO_{\eta}$, there is a constant $C = C(n, p, \eta)$ such that
$$\|[b, H_{\eta}] : L^p(\mathbb{R}^n) \to L^p(\mathbb{R}^n)\| \le C \|b\|_{BMO_{\eta}}.$$
\end{theorem} 
In a similar manner, the same testing condition used in the parabolic setting can be used here, with a slight modification in the definition to account for the fact that the cubes lie in $\mathbb{R}^n$ rather than the plane. Consider an $\eta$-cube $Q = I_1 \times \cdots \times I_n$ with side-length defined by
$$\ell(Q) = |I_1|^{\alpha_1} = \cdots = |I_n|^{\alpha_n}.$$
The testing condition is then identical to the parabolic setting, namely that
$$\|b\|_{\text{test}, \eta} := \sup_{Q \in \mathcal{P}}  \fint \left|b(x) - \frac{1}{\int_{I_{x, E_Q}} \frac{dt}{t}} \int_{I_{x, E_Q}} b(x - \eta(t))\, \frac{dt}{t}\right| \, dx < \infty$$
with $E_Q = \{x - \eta(t) : x \in Q, 9 \ell(Q) \le t \le 10 \ell(Q)\}$ and $I_{x, E_Q} = \{t : x - \eta(t) \in E_Q\}.$ On analyzing the proof of Theorem \ref{theorem:parabolic_necessary}, one finds that the only change necessary is that the equivalence $|E_Q| \sim |Q|$ of equation \eqref{eqn:similar_measures} now involves a constant depending on the curve and the dimension. This leaves us with the following:
\begin{theorem}
If the commutator $[b, H_{\eta}]$ is bounded from $L^2(\mathbb{R}^n)$ to itself, $b$ is an $\eta$-adapted testing-BMO function with
$$\|b\|_{\text{test}, \eta} \lesssim \|[b, H_{\eta}] : L^2 \to L^2\|.$$
\end{theorem}

As a final application, we may further extend the result to curves with non-vanishing torsion; that is, curves $\nu : [-1, 1] \to \mathbb{R}^n$ for which
$$\det(\nu^{(1)}(s), ..., \nu^{(n)}(s)) \ne 0$$
for all $s \in (-1, 1)$. Again, such a curve induces a family of $\nu$-cubes with the necessary properties and sparse domination theorem; in this case, however, the Hilbert transform $H_{\nu}$ must be truncated to a local operator
$$H_{\nu}^1 f(x) = p.v. \int_{-1}^1 f(x - \nu(t)) \, \frac{dt}{t}.$$ 
The details of this construction are given in \cite{ClaOu18}. We then have
\begin{theorem}
Let $p \in (1, \infty)$. If $\nu$ is a curve with non-vanishing torsion and $b \in BMO_{\nu}$, there is a constant $C = C(n, p, \nu)$ such that 
$$\|[b, H_{\nu}^1] : L^p(\mathbb{R}^n) \to L^p(\mathbb{R}^n)\| \le C \|b\|_{BMO_{\nu}}.$$
\end{theorem}

\section{Curves with torsion}

For comparison to the parabolic case, we will also consider the Hilbert transform along a line; this will show how it is important to have curvature in order to have compatibility between BMO and the commutator.

\begin{defn} The Hilbert transform along a line is defined by
$$H_{\tau} f(x, y) = \int_{-\infty}^{\infty} f(x - t, y) \, \frac{dt}{t}.$$
Alternatively, if $f^y$ denotes the section of $f$ at height $y$ (that is, the function $f^y(x) := f(x, y)$) then we also have
$$H_{\tau} f(x, y) = (Hf^y)(x)$$
where $H$ is the ordinary one-dimensional Hilbert transform.
\end{defn}

\begin{prop}
The commutator $[b, H_{\tau}]$ is bounded from $L^p(\mathbb{R}^2) \to L^p(\mathbb{R}^2)$ with norm at most
$$\|b\|_{BMO, \infty} := \|y \mapsto \|b^y\|_{BMO(dx)} \|_{\infty}.$$
\end{prop}

\begin{proof}
Following the definition of $H_{\tau}$, we can compute that
\begin{align*}
\big([b, H_{\tau}]f\big)(x, y) &= b(x, y) \big(H_{\tau}f\big)(x, y) - \big(H_{\tau}bf\big)(x, y) \\
	&= b^y(x) \big(Hf^y\big)(x) - \big(H(bf)^y\big)(x) \\
	&= \big(b^y Hf^y\big)(x) - \big(Hb^y f^y\big)(x) \\
	&= \big([b^y, H]f^y\big)(x).
\end{align*}
Using the standard one-dimensional commutator bounds,
\begin{align*}
\|[b, H_{\tau}]f\|_p^p &= \int_{-\infty}^{\infty} \int_{-\infty}^{\infty} \big|\big([b^y, H]f^y\big)(x)\big|^p \, dx \, dy \\
&= \int_{-\infty}^{\infty} \|[b^y, H] f^y\|_{L^p(\mathbb{R})} \, dy \\
&\le \int_{-\infty}^{\infty} \|b^y\|_{BMO(dx)} \|f^y\|_p^p \, dy \\
&\le \esssup_y \|b^y\|_{BMO(dx)} \|f\|_p^p
\end{align*}
as desired.
\end{proof}

Note that the estimate of this proposition is sharp because the one-dimensional bounds are sharp. We are now interested in comparing the operator norm of the commutator with the (planar) BMO norm of $b$. Unfortunately, they are not compatible.

\begin{prop}
There is a function $b \notin BMO(\mathbb{R}^2)$ with $\|b\|_{BMO, \infty} < \infty$.
\end{prop}

\begin{proof}
Consider a function $b(x, y) = g(x) + \psi(y)$ with $g \in BMO(dx)$ but $\psi \notin BMO(dx)$. Given a cube $Q = I \times J$ with $|I| = |J|$, we can compute that
$$b_Q = \frac{1}{|I||J|} \int_I \int_J g(x) \, dy \, dx + \int_J \int_I \psi(Y) \, dx \, dy = g_I + \psi_J.$$
Therefore,
\begin{align*}
\frac{1}{|Q|} \int_Q |b - b_Q| \, dA &= \frac{1}{|Q|} \int_Q |g(x) - g_I + \psi(y) - \psi_J| \, dy \, dx \\
&\ge \frac{1}{|Q|} \int_J \int_I |\psi(y) - \psi_J| \, dx \, dy - \frac{1}{|Q|} \int_I \int_J |g(x) - g_I| \, dy \, dx \\
&= \frac{1}{|J|} \int_J |\psi(y) - \psi_J| \, dy - \frac{1}{|I|} \int_I |g(x) - g_I| \, dx
\end{align*}
Now since $g \in BMO(dx)$, the second term is controlled; the first term can be made arbitrarily large, however, and the result follows.
\end{proof}

The above proposition shows that our new space is not necessarily a subspace of planar BMO; the reverse inclusion does not hold either.\begin{prop}
There is a function $b \in BMO(\mathbb{R}^2)$ such that $\|b\|_{BMO, \infty} = \infty$.
\end{prop}

\begin{proof}
Consider a function $b(x, y) = g(x) \psi(y)$; we can easily compute
$$b_Q = g_I \psi_J$$
in a manner similar to the previous proposition. We may also compute
$$\|b^y\|_{BMO(dx)} = \|\psi(y) g(x)\|_{BMO(dx)} = |\psi(y)| \|g\|_{BMO(dx)}$$
from which it follows that
$$\|b\|_{BMO, \infty} = \|\psi\|_{\infty} \|g\|_{BMO(dx)}.$$
However, we now have
\begin{align*}
\frac{1}{|Q|} \int_Q |b - b_Q| \, dA &= \frac{1}{|I||J|} \int_J \int_I |g(x) \psi(y) - g_I \psi_J| \, dx \, dy \\
&= \frac{1}{|I||J|} \int_J \int_I |g(x) \psi(y) - g(x) \psi_J + g(x) \psi_J - g_I \psi_J| \, dx \, dy \\
&\le \frac{1}{|I||J|} \int_J \int_I |g(x)| |\psi(y) - \psi_J| \, dx \, dy + \\
&\quad + \frac{1}{|I||J|} \int_J \int_I |\psi_J| |g(x) - g_I| \, dx \, dy \\
&= \frac{1}{|I|} \int_I |g| \cdot \frac{1}{|J|} \int_J |\psi - \psi_J|  + \frac{1}{|J|} \int_J |\psi| \cdot \frac{1}{|I|} \int_I |g - g_I|
\end{align*}
Now if $\psi$ is unbounded but $g \in L^{\infty}$ is non-constant, we have $\|b\|_{BMO, \infty} = \infty$. If $\psi \in BMO$, then 
$$ \frac{1}{|I|} \int_I |g| \cdot \frac{1}{|J|} \int_J |\psi - \psi_J| $$
is uniformly bounded over all $I$ and $J$. Finally, choose $\psi$ so that the averages of $\psi$ blow up slowly as $|J| \to 0$, but $g$ so that the oscillation of $g$ over $I$ tends to zero more quickly as $|I| \to 0$ (for example, this can be done with $\psi \in LMO$ and $g \in L^{\infty} \subseteq VMO$); this gives $\|b\|_{BMO} < \infty$. 
\end{proof}

\bigskip

{\bf Acknowledgement:} Z. Guo and J. Li are supported by ARC DP 170101060. J. Li is also supported by
Macquarie University Research Seeding Grant. 
 B. D. Wick's research is partially supported by National Science Foundation -- DMS \# 1800057 and DMS \# 1560955 and Australian Research Council -- DP 190100970.

\section{References}

\begingroup
\renewcommand{\section}[2]{}%
\bibliographystyle{plain}
\bibliography{commutator_bibliography}
\endgroup

\end{document}